\documentclass[a4paper,11pt]{amsart}

\usepackage{amscd,amsfonts,amsmath,amssymb}

\input xy
\xyoption{all}
\usepackage[all]{xy}

\newcommand{\cal}{\mathcal}

\newcommand{\C}{{\mathbb C}}
\newcommand{\Q}{{\mathbb Q}}

\newcommand{\lb}{\lbrace}
\newcommand{\rb}{\rbrace}

\newcommand{\Z}{\ensuremath{\mathbb{Z}}}

\newcounter{subsub}[subsection]

\def\hfl#1#2{\smash{\mathop{\hbox to 12mm{\rightarrowfill}}
\limits^{\scriptstyle#1}_{\scriptstyle#2}}}

\long\def\InsertFig#1 #2 #3 #4\EndFig{
\hbox{\hskip #1 mm$\vbox to #2 mm{\vfil\includegraphics{#3}}#4$}}
\long\def\LabelTeX#1 #2 #3\ELTX{\rlap{\kern#1mm\raise#2mm\hbox{#3}}}

\newtheorem{theorem}{Th\'eor\`eme}[section]
\newtheorem{lemma}[theorem]{Lemme}
\newtheorem{example}[theorem]{Exemple}
\newtheorem{corollary}[theorem]{Corollaire}
\newtheorem{proposition}[theorem]{Proposition}
\newtheorem{question}[theorem]{Question}
\newtheorem{remark}[theorem]{Remarque}

\newtheorem{definition}[theorem]{D\'efinition}
\newtheorem{conjecture}[theorem]{Conjecture}

\newtheorem{abundance conjecture}[theorem]{Abundance Conjecture}
\newtheorem{simplified lemma}[theorem]{Simplified Lemma}

\newtheorem{main questions}[theorem]{Main Questions}

\newcommand\ra{{\rightarrow}}

\begin{document}
\title{ORBIFOLDES \`A PREMI\`ERE CLASSE DE CHERN NULLE}
\author{Fr\'ed\'eric CAMPANA}
%\date{2002}
\maketitle
{\medskip}

\section{\bf INTRODUCTION.}

On \'etablit ici une version orbifolde du th\'eor\`eme de d\'ecomposition de Bogomolov 
pour les vari\'et\'es K\" ahl\'eriennes compactes 
\`a premi\`ere classe de Chern nulle (voir [Bo74], [Bo78], [Y78], [Bea83] et [P97] pour certains cas singuliers). Les 
orbifoldes K\" ahl\'eriennes consid\'er\'ees sont les V-vari\'et\'es de Satake, ou encore les 
vari\'et\'es \`a singularit\'es quotient. 

La d\'emonstration (voir \S 2, 3, 4) est une extension orbifolde directe du cas lisse, qui repose sur 
la construction de m\'etriques de K\" ahler Ricci-plates, le scindage de Cheeger-Gromoll, et le th\'eor\`eme 
de d\'ecomposition de De Rham. 

Appliqu\'e au cas des surfaces $K3$ normales, on en d\'eduit en particulier le fait que le 
groupe fondamental du lieu lisse d'une telle surface est soit fini, soit extension fini de $\Bbb Z^{\oplus 4}$ par un groupe fini. 
Cette propri\'et\'e avait \'et\'e conjectur\'ee dans [K-Z00]. Voir la bibliographie de cet article, ainsi que [KZ02], [S-Z01], et [C-K-O 03]
pour de nombreux cas \'etablis par voie alg\'ebro-g\'eom\'etrique.

 Notre motivation \'etait de v\'erifier 
la conjecture formul\'ee dans [C01] de la presque-ab\'elianit\'e du groupe fondamental et l'annulation de la pseudo-m\'etrique de Kobayashi 
pour 
certaines vari\'et\'es K\" ahl\'eriennes sp\'eciales de dimension trois. Voir \S 7 ci-dessous pour plus de d\'etails.

On pourra trouver dans [P97] et [A03] des r\'esultats \'etablissant une telle d\'ecomposition dans des situations diff\'erentes. ([A] traite par voie alg\'ebro-g\'eom\'etrique le cas du morphisme d'Albanese pour une paire log-terminale projective $(X,B)$).

Je voudrais remercier J.P. Demailly pour ses indications sur la preuve du th\'eor\`eme \ref{ke}, ainsi que J. Maubon, 
qui m'a fourni les r\'ef\'erences [B93] et [B-Z94] sur la version orbifolde du th\'eor\`eme 
de scindage de Cheeger-Gromoll, qui joue un r\^ole crucial ici.
\bigskip

\section{\bf SUMMARY.}

We establish an orbifold version of Bogomolov decomposition theorem for compact K\" ahler manifolds with trivial first Chern class. The 
orbifolds here considered are Satake's V-manifolds. The proof is a direct orbifold extension of the proof in the smooth case, obtained 
by using Ricci-flat K\" ahler metrics, Cheeger-Gromoll splitting theorem, and De Rham decomposition theorem.

Applied to normal $K3$ surfaces, we thus obtain in particular the fact that the fundamental group of their smooth locus is either finite, or an 
extension of $\Bbb Z^{\oplus 4}$ by a finite group, in which case the $K3$ under consideration is uniformised (in the orbifold sense), by a 
complex torus (of dimension $2$). This property was conjectured in [K-Z00, 3.12], and proved in many situations by algebro-geometric methods. 
We refer to [K-Z02], [S-Z01], and [C-K-O 03] for further details on the known cases.

Our motivation was to check conjectures of almost abelianity and vanishing of 
the Kobayashi pseudo-metric for certain special threefolds. See \S 7 for further details.

\section{\bf ORBIFOLDES PURES.}

Nous suivons ici (\`a la terminologie pr\`es) les d\'efinitions de [D-K00, \S 6]. Les notions introduites ici
apparaissent dans [S56] auquel nous renvoyons pour plus de d\'etails. 

\begin{definition} Une {\bf orbifolde pure} est un espace analytique complexe normal $X$ n'ayant que des singularit\'es 
quotient. 
\end{definition}
Pour chaque point $a\in X$, on a donc un voisinage $U$ de $a$ dans $X$, et un isomorphisme 
$\phi:U\ra \tilde{U}/G$, dans lequel $\tilde{U}$ est un ouvert de $\C^n$, et $G\subset Gl(n,\C)$ un sous-groupe fini 
agissant sur $\tilde{U}$, avec $a=\phi (0)$. On peut choisir $G$ tel que les points fixes de tout $g\in G$ ait 
un ensemble de points fixes de codimension au moins $2$. Un tel triplet $(U,\tilde{U}/G,\phi)$ sera appel\'e une {\bf uniformisation 
locale} de $X$ en $a$. Remarquons que, si l'on restreint l'ouvert $U$, le groupe $G$ est alors le groupe fondamental local 
$\pi_1^{loc}(a,X)$ de $X$ en $a$, et que l'uniformisation pr\'ec\'edente est alors unique \`a unique isomorphisme pr\`es.

On notera $X^*\subset X$ (resp. $Sing(X)\subset X$) le lieu lisse (resp. singulier) de $X$.

Les orbifoldes pures sont simplement appel\'ees orbifoldes dans [D-K00]. On ajoute ici l'adjectif ``pur" pour les 
distinguer des orbifoldes plus g\'en\'erales de [C01], dans lesquelles la structure d'orbifolde comporte \'egalement 
un $\Q$-diviseur $\Delta$.

Les objets $\Theta$ sur l'orbifolde pure $X$ (tels que le faisceau structural, le faisceau canonique, les m\'etriques de K\" ahler, 
la forme de courbure d'une telle m\'etrique) sont les objets usuels sur la partie lisse de $X$, avec la condition que, dans toute 
uniformisation locale $(U,\tilde{U}/G,\phi)$ de $X$ en $a$, l'objet $\Theta$ soit sur la partie lisse de $U$, l'image par 
$\phi$ d'un objet correspondant $G$-invariant $\tilde{\Theta}$ sur $\tilde{U}$. Si $\alpha$ est une forme de degr\'e maximum sur 
l'ouvert $U$ d'une uniformisation locale, on d\'efinit ainsi: $\int_{U}\alpha:=(1/Card(G)).\int_{\tilde{U}}\tilde{\alpha}$.

Par exemple: $\cal O_U:=p_*(\cal O_V)$, et $ K_U^m:=p_*( K_V^m),\forall m\in \Z$. Ce dernier faisceau est donc localement libre 
si $m$ est un multiple de $Card(G)$. Donc $K_X^m$ est localement libre si $X$ est compacte et $m$ suffisamment divisible.

On peut donc d\'efinir la premi\`ere classe de Chern $c_1(X):=-(1/m).(c_1(K_X^m))$, si $X$ est compacte, et $m$ comme ci-dessus. 
Cette classe de Chern peut encore, comme dans le cas lisse, 
\^etre calcul\'ee comme la classe de cohomologie de la forme de courbure de la connexion 
d'une m\'etrique hermitienne sur le fibr\'e $K_X^{-m}$, pour $m$ suffisamment divisible.

\begin{example} Soit $X$ une orbifolde pure compacte, et de dimension $2$. Alors $Sing(X)\subset X$ est un ensemble fini. Soit 
$r:X'\ra X$ une r\'esolution de $X$. Si $X'$ est K\" ahler, alors $X$ est K\" ahler. 
En effet, la m\'etrique de $X'$ d\'efinit sur 
$X^*$ une m\'etrique de K\" ahler que l'on peut recoller 
avec une m\'etrique de K\" ahler locale d\'efinie sur une uniformisation locale 
de chacun des points de $Sing(X)$.
 
Soit $X'$ une surface $K3$ et $D$ un diviseur r\'eduit effectif exceptionnel de $X'$. Soit $r:X'\ra X$ la contraction de 
$D$ dans $X'$. Alors (par [K-M98; 4.20], [B-P-V84; II.7.3, III.2.4], et [D79]),
 $X$ est une orbifolde pure compacte connexe K\" ahler de dimension $2$ avec fibr\'e canonique trivial (et $c_1(X)=0$).
R\'eciproquement, si $X'$ est une orbifolde pure singuli\`ere 
compacte connexe et K\" ahler de dimension $2$ \`a fibr\'e canonique trivial, on 
d\'eduit de [K-M98, 4.20] \`a nouveau que $X'$, dans  la r\'esolution minimale $r:X'\ra X$, est une surface $K3$.
\end{example}

Si $X$ est une orbifolde pure, munie d'une m\'etrique orbifolde de  K\" ahler, on a encore l'\'egalit\'e 
$2.\Delta'=2.\Delta"=\Delta$, pour les Laplaciens (orbifoldes) d\'eduits de cette m\'etrique, ceci par 
les m\^emes arguments (locaux) que dans le cas lisse. En particulier, le lemme du 
$\partial\overline{\partial}$ est encore vrai, dans 
le contexte orbifolde. Une $(1,1)$-forme $\alpha$ r\'eelle $d$-exacte est donc de la forme: 
$\alpha=i.\partial\overline{\partial}f$, pour $f:X\ra \Bbb R$ de classe $\cal C^{\infty}$ (au sens orbifolde).

\bigskip
\section{\bf M\'ETRIQUES DE K\" AHLER-EINSTEIN.}

Lorsque $X$ est une orbifolde compacte avec $K_X$ ample (ou $c_1(X)<0$), il est \'etabli dans [K87] que $X$ 
admet une m\'etrique de K\" ahler-Einstein. La d\'emonstration est la m\^eme que celle ([Au78]) du cas lisse. Le cas des 
orbifoldes de Fano ($K_X^{-1}$ ample, ou $c_1(X)>0$) est \'etudi\'e dans [D-K00]. On va ici constater que le cas o\`u $c_1=0$ 
est, dans le contexte orbifolde, exactement similaire \`a celui du cas lisse. (Curieusement, il nous a \'et\'e 
impossible de trouver cette observation dans la litt\'erature).

\begin{theorem}\label{ke} Soit $X$ une orbifolde pure K\"ahl\'erienne compacte et connexe avec $c_1(X)=0$. 
Toute classe de K\" ahler de $X$ est repr\'esent\'ee par
une unique m\'etrique de K\" ahler Ricci-plate sur $X$.
\end{theorem}
{\bf D\'emonstration:} C'est une adaptation imm\'ediate de la d\'emonstration dans le cas lisse, 
telle qu'elle est expos\'ee, par exemple, dans [Siu 87; pp. 85-113]. Dans une premi\`ere \'etape, on se ram\`ene \`a 
la r\'esolution de l'\'equation de Monge-Amp\`ere: $det(g_{i\bar j}+\phi_{i\bar j})=e^{F}.det(g_{i\bar j})$, dans laquelle 
$(g_{i\bar j})$ est (en coordonn\'ees locales) une m\'etrique de K\"ahler de la classe donn\'ee d'une forme de K\" ahler $\omega$,
 $\phi:X\ra \Bbb R$ est une fonction 
$\cal C^{\infty}$ inconnue, $\phi_{i\bar j}:=i\partial_i \overline{\partial}Ñ{\bar j}\phi$, et: 
$i\partial \overline {\partial}F=R_{i\bar j}=-i\partial_i\overline{\partial}_{\bar j}log (det ((g_{l\bar k}))$, puisque la 
premi\`ere clase de Chern de $X$ est nulle. 
Le lemme du $\partial \overline{\partial}$ orbifolde fournit cette mise en \'equation.

On applique alors la m\'ethode de continuit\'e, consid\'erant la famille d' \'equations $(*)_t$:
 
$det(g_{i\bar j}+\phi_{i\bar j})=A_t.e^{t.F}.det(g_{i\bar j})$, $t\in [0,1]$, $A_t:=(vol_g(M)). (\int _M e^{tF})^{-1}$. Pour $t=0$, 
on a une solution $\phi=0$, et la m\'etrique $g$. Pour $t=1$, on a donc une solution du probl\`eme. 
On consid\`ere l'ensemble non vide des $t$ pour lesquels l'\'equation $(*)_t$ a une solution. 

Cet ensemble est ouvert, par une application du th\'eor\`eme des fonctions implicites et lin\'earisation du probl\`eme. 
La d\'emonstration est celle du cas lisse, en regardant les solutions localement dans les ouverts $\tilde U$. On peut 
consulter, par exemple, [D-K00, \S 6.2] pour le cas analogue des orbifoldes de Fano.

La partie la plus d\'elicate est de d\'emontrer que cet ensemble est ferm\'e. Il suffit pour cela d'\'etablir des estimations
 a priori pour les solutions $\phi_t$  de $(*)_t$, portant sur les normes uniformes de $\phi_t$, puis 
sur celles des d\'eriv\'ees du second ordre, puis sur les normes de H\" older $C^{2,\alpha}$ de $\phi_t$.

Les estim\'ees uniformes sont obtenues dans [Siu87;\S 2] par la m\'ethode d'it\'eration de Moser. La seule chose \`a v\'erifier
 par rapport au cas lisse est la validit\'e de la formule d'int\'egration par parties dans le cas orbifolde, ce qui est imm\'ediat. 

Les estim\'ees suivantes, (bas\'ees sur le principe du maximum et une in\'egalit\'e de Harnack due \`a Moser) 
qui sont les plus difficiles (voir [Siu 87, \S 4 et \S 5]), r\'esultent d'arguments locaux, et peuvent donc
\^etre trait\'ees dans des ouverts uniformisants $\tilde U$, sans changement. $\square$

\begin{corollary} Soit $r:X'\ra X$ la contraction d'un diviseur exceptionnel $D'\subset X'$ d'une surface $K3$ $X'$. Alors 
$X$ est une orbifolde K\" ahl\'erienne pure \`a fibr\'e canonique trivial. 
Toute classe de K\" ahler sur $X$ est repr\'esent\'ee par une 
unique m\'etrique de K\" ahler Ricci-plate.
\end{corollary}

\bigskip

\section{\bf UNIFORMISATION ORBIFOLDE.}

On va consid\'erer ici une orbifolde pure connexe $X$, dont on notera $X^*$ le lieu lisse. On notera $\pi_1^{orb}(X):=\pi_1(X^*)$ 
son {\bf groupe fondamental d'orbifolde}. On a une surjection naturelle $\pi_1^{orb}(X)\ra \pi_1(X)$, qui n'est en 
g\'en\'eral pas injective (consid\'erer par example une surface $K3$ normale de Kummer $X$, 
quotient d'un tore complexe de dimension de dimension $2$. Alors $X$ est simplement connexe, mais son 
groupe fondamental d'orbifolde est 
isomorphe \`a $\Bbb Z^{\oplus 4}$). 

Nous donnons quelques d\'etails sur ces notions faciles, faute de r\'ef\'erences 
accessibles sur ce sujet (les r\'ef\'erences donn\'ees dans [B-Z94] renvoient toutes \`a des  
th\`eses, ou publications internes \`a des universit\'es).

On va \'etablir une relation entre le groupe fondamental d'orbifolde de $X$ et les rev\^etements orbifoldes de 
$X$, similaire \`a celle du cas lisse.

\begin{definition}\label{drevorb} Un rev\^etement d'orbifolde $r:X'\ra X$ de $X$ est la donn\'ee d'une orbifolde pure et connexe 
$X'$ et d'une application holomorphe $r$ telle que:

{\bf (1)} Si $X'':=r^{-1}(X^*)$, 
et si $r^*:X''\ra X^*$ est la restriction de $r$ \`a $X''$, alors $r^*$ est un rev\^etement \'etale de $X^*$.

{\bf (2)} Pour tout $a\in X$, muni d'une uniformisation locale $(U,\tilde{U}/G,\phi)$, avec $\tilde{U}$ simplement connexe, 
et pour toute composante connexe $U"$ de $(r^{*})^{-1}(U\cap X^*)$, il existe $a'\in r^{-1}(a)$, un sous-groupe $G'\subset G$, 
et une uniformisation locale 
$(U', \tilde{U}/G', \phi')$ de $X'$ en $a'$, telle que $U"=U'\cap X"$, et telle que:
$\phi^{-1}\circ r\circ \phi':\tilde{U}/G'\ra \tilde{U}/G$ soit le quotient naturel induit par l'inclusion de $G'$ dans $G$.
\end{definition}

\
On d\'efinit comme dans le cas lisse la notion d'isomorphisme de rev\^{etements orbifoldes.}¥

Tout comme dans le cas lisse, on a une bijection naturelle entre sous-groupes 
du groupe fondamental orbifolde et (classes d'isomorphismes de)¥rev\^etements orbifoldes. Plus pr\'ecis\'ement:

\begin{proposition} Soit $X$ une orbifolde pure et connexe, et $b\in X^*$ un point-base 
lisse de $X$. L'application naturelle entre les classes d'isomorphisme de rev\^etements orbifoldes connexes de $X$ et les sous-groupes 
de $\pi_1^{orb}(X,b):=\pi_1(X^*,b)$, qui associe \`a un rev\^etement orbifolde 
le sous-groupe d\'efini par sa restriction au-dessus de $X^*$ (d\'efinie par \ref{drevorb} {\bf (1)} ci-dessus) est bijective. 
\end{proposition}

{\bf D\'emonstration:} L'injectivit\'e est \'evidente. La surjectivit\'e provient de ce que l'on peut partiellement 
compactifier tout rev\^etement connexe $r^*:X"\ra X^*$ de $X^*$ en un rev\^etement orbifolde $r:X'\ra X$ 
de $X$, en utilisant la finitude 
des groupes fondamentaux locaux $\pi_1^{loc}(a,X)$. 
Plus pr\'ecis\'ement, si $a\in Sing(X)$, et si $(U, \tilde{U}/G,\phi)$ est une uniformisation 
locale de $X$ en $a$ telle que $\tilde{U}$ soit simplement connexe, alors chaque composante connexe de 
$(r^*)^{-1}(U^*)$ est isomorphe 
\`a $(\tilde{U}/G'):=U'$, pour un unique sous-groupe $G'\subset G$ (un point-base de $U$ ayant \'et\'e fix\'e, notant: 
$U^*:=(X^*\cap U)$). 
On peut alors recoller $(r^*)^{-1}(U^*)$ et $\Gamma\times U'$, gr\^ace aux isomorphismes pr\'ec\'edents, 
si $\Gamma$ est l'ensemble des composantes connexes de $(r^*)^{-1}(U^*)$. 
Les d\'etails ne pr\'esentent aucune difficult\'e $\square$

\begin{definition} Le rev\^etement $\bar{r}:\bar{X}\ra X$ qui correspond au sous-groupe trivial $\lb 1\rb$ de $\pi_1^{ob}(X)$ 
dans la bijection pr\'ec\'edente est appel\'e le {\bf rev\^etement universel orbifolde } de $X$. 
\end{definition}

On remarquera que $\bar{X}$ 
 est lisse (et fournit donc une uniformisation globale de $X$) si et seulement si, pour chaque $a\in X$, le morphisme naturel 
$\pi_1^{loc}(a,X)\ra \pi_1^{orb}(X)$ d\'eduit 
de l'injection de $U^*$ dans $X^*$, est injectif. (Je remercie T. Delzant pour m'avoir signal\'e ce point).

Remarquons aussi que, si $X$ est une orbifolde compacte munie d'une m\'etrique d'orbifolde Riemannienne arbitraire, 
alors: $\bar{X}$, munie de la m\'etrique image r\'eciproque de celle de $X$ par $\bar{r}$, est compl\`ete. 

De plus, si $\bar{X}$ est non-compacte, elle contient (au moins) une droite (ie: une g\'eod\'esique 
d\'efinie sur $\Bbb R$, minimisant la distance entre deux quelconques de ses points). Voir [B-Z94] pour cette assertion.

Nous allons maintenant fournir la version orbifolde du th\'eor\`eme de d\'ecomposition 
Riemannien de De Rham. Ce r\'esultat 
est utilis\'e sans d\'emonstration dans [C-P02]. 

Nous \'etablissons ici l' observation qui permet d'adapter la 
d\'emonstration du cas lisse. Pour les r\'esultats classiques relatifs aux groupes d'holonomie, 
nous renvoyons \`a [Be87; \S 10], et \`a sa bibliographie.

\begin{proposition}\label{tdecorb} Soit $X$ une orbifolde pure munie d'une m\'etrique Riemannienne compl\`ete. 
Soit $X^*$ sa partie lisse, suppos\'ee simplement connexe. Alors $X$ admet une 
d\'ecomposition en produit Riemannien fini d'orbifoldes $X=\Pi_{j\in J} X_j$ tel que: $X^*=\Pi _{j\in J} X_i^*$, l'action du 
groupe d'holonomie de chacun des facteurs $X_j^*$ \'etant irr\'eductible.
\end{proposition}

{\bf D\'emonstration:} Le groupe d'holonomie de $X^*$ (relatif \`a un point fix\'e $b\in X^*$) est un 
sous-groupe de Lie du groupe 
orthogonal appropri\'e. En particulier, il y est ferm\'e. Ce groupe d'holonomie, \'etant ferm\'e et engendr\'e 
par les d\'eplacements parall\`eles le long de ``petits lassos" (au sens de Kobayashi-Nomizu) contenus dans $X^*$, contient donc 
(\`a conjugaison pr\`es) le groupe d'holonomie de $\tilde{U}$, si $(U,\tilde{U}/G,\phi)$ est une uniformisation locale 
de $X$ en un point arbitraire $a\in X$ (ceci parce que ce groupe est engendr\'e par les 
d\'eplacements parall\`eles le long des ``petits lassos" contenus dans $U$). 

Or le d\'eplacement parall\`ele $T_c$ le long d'un lacet $c$ de $U$ 
se rel\`eve en le d\'eplacement parall\`ele $T_{c'}$ le long du lacet
$c':=u"\circ \tilde{c}\circ u'$, dans lequel $u"$ (resp. $u'$) 
est le rel\`evement \`a $\tilde{U}$ de $u^{-1}$ (resp. $u$) partant de $0\in \tilde{U}$ et aboutissant \`a 
$d'$ (resp. $d"$), point initial (resp. final) de $\tilde{c}$, o\`u $u$  est un chemin dans $U$ qui joint le point $a$ 
au point  initial $d$ de $c$, et o\`u $\tilde{c}$ rel\`eve $c$ \`a $\tilde{U}$.

Cette observation montre que la d\'ecomposition locale, $ G$-invariante, en $0\in\tilde{U}$ de $\tilde{U}$ en produit 
de facteurs irr\'eductibles pour l'holonomie de $\tilde{U}$ est compatible avec les 
(ie: fournit des sous-facteurs des) feuilletages 
d\'efinissant les facteurs irr\'eductibles pour l'holonomie de $X^*$. 

A partir de ce point, la d\'emonstration de la version orbifolde du th\'eor\`eme de De Rham utilise la compl\'etion de $X$ de 
 la m\^eme mani\`ere que dans le cas lisse
$\square$

\bigskip

\section{\bf LE TH\'EOR\`EME DE SCINDAGE ORBIFOLDE.}

Nous avons la version suivante du th\'eor\`eme de scindage de Cheeger-Gromoll. Je remercie vivement J.Maubon qui m'a appris 
l'existence de cette r\'ef\'erence.

\begin{theorem} ([Bo-Z94] Soit $\tilde{X}$ une orbifolde Riemannienne compl\`ete 
\`a courbure de Ricci partout positive ou nulle. Alors: 
$\tilde{X}=\tilde{N}\times \Bbb R^l$, o\`u $N$ est une orbifolde Riemannienne compl\`ete ne contenant pas de 
droite. (Le produit pr\'ec\'edent est pris au sens Riemannien, 
la m\'etrique sur $\Bbb R^l$ \'etant la m\'etrique plate).
\end{theorem}

Remarquons que la notion d'orbifolde utilis\'ee dans [B-Z94] (et que nous ne d\'efinissons pas ici) 
est plus g\'en\'erale que la notion d'orbifolde pure consid\'er\'ee ci-dessus. Ce r\'esultat s'applique donc en particulier aux 
orbifoldes K\" ahl\'eriennes Ricci-plates de la \S 3 ci-dessus, et fournit donc le: 

\begin{corollary} ( [Bo-Z94]) Soit $X$ une orbifolde Riemannienne compacte \`a courbure de Ricci 
partout positive ou nulle. Alors $X$ a 
un rev\^etement d'orbifolde fini $\bar X$ qui est un produit Riemannien $\bar X=\bar N\times T$, dans lequel $\bar N$ est une 
orbifolde simplement connexe (au sens orbifolde), et $T$ est un tore (r\'eel) plat.

En particulier, $\pi_1^{orb}(X):=G$ est extension de $\Bbb Z^{\oplus l}$ par un sous-groupe fini $F$ de $G$ (On dira que $G$ est 
{\bf presque-ab\'elien de rang $l$}).
\end{corollary}

La formulation de ce r\'esultat n'est pas donn\'ee dans [B-Z94], mais s'en d\'eduit imm\'ediatement \`a l'aide du th\'eor\'eme 
de Bieberbach (un groupe cristallographique a un sous-groupe d'indice fini sans point fixe, donc ab\'elien).

\begin{corollary}\label{cuniforb} Soit $X$ une orbifolde K\"ahl\'erienne compacte et connexe avec $c_1(X)=0$. Alors: 

{\bf (1)} Le groupe fondamental de $X$ est presque ab\'elien de rang $2m$, $m\leq dim(X)$.

{\bf (2)} $X$ admet un rev\^etement d'orbifolde fini $\bar X=\bar N\times T$, dans lequel $\bar N$ est une 
orbifolde K\" ahl\'erienne compacte et simplement connexe au sens orbifolde, et $T$ est un tore complexe de dimension $m$.
\end{corollary}

Ce corollaire est un cas particulier du corollaire \ref{cuniforb} ci-dessous. Mais on peut le d\'eduire plus \'el\'ementairement de 
l'observation initiale de la d\'emonstration de \ref{tdecorb} ci-dessus, et du th\'eor\`eme de Lichn\'erowicz ([Li55, p. 264]), 
qui affirme que les facteurs locaux de la d\'ecomposition irr\'eductibe d'holonomie d'une vari\'et\'es K\" ahl\'erienne 
sont encore K\" ahl\'eriens (et en particulier complexes).

\medskip

Le corollaire \ref{cuniforb} {\bf (2)} ci-dessus peut \^etre consid\'erablement pr\'ecis\'e par \ref{t} ci-dessous, dont la 
d\'emonstration consiste simplement 
\`a appliquer les arguments de [Bea83], utilisant la version orbifolde \ref{tdecorb} du th\'eor\`eme de d\'ecomposition de De Rham.
Voir [P97] pour des r\'esultats de d\'ecomposition similaires pour les vari\'et\'es projectives $\Bbb Q$-factorielles 
\`a singularit\'es terminales.
 
\begin{theorem}\label{t} Soit $X$ une orbifolde K\"ahl\'erienne compacte et connexe avec $c_1(X)=0$. Alors: 
$X$ admet un rev\^etement d'orbifolde fini $\bar X=\bar C\times \bar S\times T$, dans lequel $\bar C$ (resp. $\bar S$) est un 
produit fini d'orbifoldes K\"ahl\'eriennes de {\bf Calabi-Yau} (resp. {\bf Hyperk\" ahl\'eriennes}), et $T$ est un tore 
complexe de dimension $m$. Cette d\'ecomposition est unique.
\end{theorem}

On d\'efinit les notions qui apparaissent dans le r\'esultat pr\'ec\'edent:

\begin{definition} Une orbifolde pure $X$ de dimension complexe $m$ est dite de {\bf Calabi-Yau} (resp. {\bf Hyperk\" al\'erienne}) 
si elle est compacte, simplement connexe au sens orbifolde, et si elle admet une m\'etrique K\"ahl\'erienne Ricci plate 
dont la restriction \`a sa partie lisse est d'holonomie irr\'eductible \'egale \`a $SU(m)$ (resp. $Sp(m/2)$).
\end{definition}

Les arguments de [Bea83], bas\'es sur le th\'eor\`eme de Calabi-Yau, le principe d'holonomie, et la d\'ecomposition de De Rham, 
s'\'etendent donc au cas orbifolde, par ce qui pr\'ec\`ede, et fournissent:

\begin{proposition} Soit $X$ une orbifolde pure et K\" ahl\'erienne compacte connexe et simplement connexe au sens 
orbifolde. On a \'equivalence entre les propri\'et\'es {\bf (1)},{\bf (2)} suivantes:

{\bf (1)} $X$ est de Calabi-Yau (resp. Hyperk\" ahl\'erienne).

{\bf (2)} $K_X$ est trivial et $H^0(X,\Omega_X^2)=0$ (resp. $H^0(X,\Omega_X^2)$ est engendr\'e par une section partout de rang maximum sur $X^*$).
\end{proposition}

Voir aussi [Fu83] pour d'autres r\'esultats et exemples concernant les orbifoldes pures hyperk\" al\'eriennes, 
et leurs d\'eformations.

\medskip

Dans le cas des surfaces, on peut \^etre plus pr\'ecis:

\begin{corollary}\label{unifk3} Soit $r:X'\ra X$ la contraction d'un diviseur exceptionnel $D'\subset X'$, 
o\`u $X'$ est une surface $K3$. Alors le groupe fondamental d'orbifolde de $X$ est soit fini, soit presque-ab\'elien 
de rang $4$. De plus, $X$ admet un rev\^etement d'orbifolde fini qui est soit une surface 
$K3$ normale et simplement connexe (au sens orbifolde) dans le premier cas, soit un tore complexe 
de dimension $2$ dans le second cas.
\end{corollary}

\begin{remark} Le corollaire pr\'ec\'edent r\'esoud donc la conjecture (3.12) de [K-Z00] (\`a l'exception du fait qu'une 
orbifolde $K3$ simplement connexe au sens orbifolde y est aussi conjectur\'ee avoir des singularit\'es de Du Val). 

La motivation de ce r\'esultat provient ici  
des conjectures concernant les vari\'et\'es sp\'eciales formul\'ees dans [C01]. Voir \S 7. ci-dessous.

Remarquons aussi qu'il existe des surfaces $K3$ normales et simplement connexes au sens des orbifoldes, mais non lisses. 
Par example, une surface $K3$ (lisse) de Kummer dans laquelle on contracte l'une des $16$ courbes exceptionnelles. 

On dira que l'orbifolde $X$ est {\bf bonne} si son rev\^etement universel d'orbifolde est lisse. L'exemple pr\'ec\'edent est 
donc une orbifolde $K3$ qui n'est pas bonne.
\end{remark}

\begin{corollary} \label{k3orb} Soit $X'$ une surface $K3$ lisse, et $D'\subset X'$ un diviseur exceptionnel de $X'$. Alors:

{\bf (1)} Il existe une modification $T'$ d'un tore complexe ou d'une surface $K3$ lisse
$T$ de dimension $2$ et une application holomorphe 
$r:T'\ra X'$ g\'en\'eriquement finie et finie et non-ramifi\'ee au-dessus 
de $X^*:=X'-D'$,

{\bf (2)} Si l'orbifolde $X$ obtenue par contraction de $D'$ dans $X'$ est bonne et projective, alors $X'$ contient deux 
familles \`a un param\`etre de courbes 
elliptiques (singuli\`eres, en g\'en\'eral) dont le membre g\'en\'erique est disjoint de $D'$.
\end{corollary}

{\bf D\'emonstration:} La premi\`ere assertion est simplement une reformulation de \ref{unifk3}. La seconde est obtenue 
en prenant l'image par $r$ des familles \`a un param\`etre de courbes elliptiques sur $T'=T$ (qui est alors 
une surface $K3$ lisse) construites dans [Mo-Mu82], et en utilisant le 
lemme suivant, affirmant que de telles familles de courbes elliptiques n'ont pas de point base.
$\square$

\begin{lemma} Soit $X$ une surface $K3$ lisse, et $(E_s)_{s\in S}$ une famille \`a un param\`etre de courbes elliptiques de $X$. 
Alors, par chaque point de $X$, il ne passe qu'un nombre fini de courbes $E_s$.
\end{lemma}

{\bf D\'emonstration:} Soit $w:Z\ra S$ une surface elliptique obtenue par 
d\'esingularisation du graphe d'incidence de la famille $(E_s)_{s\in S}$. Soit $h:Z\ra X$ l'application holomorphe 
obtenue par composition avec la projection sur $X$ du graphe d'incidence de cette famille. Soit $v$ un g\'en\'erateur 
de l'espace vectoriel complexe $H^0(X,\Omega_X^2)$, et $v':=h^*(v)$. Alors $v'$ est une section de $K_Z$, et son lieu des z\'eros 
est donc contenu dans une r\'eunion (finie) de fibres de $w$. Si $a\in X$ \'etait un point base de la famille, $v'$ s'annulerait 
sur la multisection $h^{-1}(a)$ de $w$. Contradiction $\square$

\begin{question} L'assertion (2) de \ref{k3orb} ci-dessus reste-t'elle vraie si $X$ n'est pas suppos\'ee bonne?
\end{question}

\bigskip

\section{\bf VARI\'ET\'ES SP\'ECIALES DE DIMENSION TROIS.}

La motivation initiale des r\'esultats pr\'ec\'edents est la d\'emonstration de cas 
particuliers en dimension $3$ de conjectures formul\'ees dans [C01] concernant les vari\'et\'es sp\'eciales (voir [C01] 
pour la d\'efinition).

\begin{conjecture} Soit $Y$ une vari\'et\'e K\" ahl\'erienne compacte sp\'eciale. Aors:

{\bf (1)} $\pi_1(Y)$ est presque ab\'elien (ie: admet un sous-groupe d'indice fini ab\'elien)

{\bf (2)} La pseudo-m\'etrique de Kobayashi $d_Y$ de $Y$ est identiquement nulle sur $Y$ (ie: 

$d_Y(y,y')=0, \forall y,y'\in Y$).
\end{conjecture}

Lorsque $Y$ est de dimension $3$, la conjecture ${\bf (1)}$ est d\'emontr\'ee dans [C01], sauf dans le cas o\`u $\kappa(Y)=2$, 
et o\`u le "coeur" $c_Y:Y\ra C(Y)$ de $Y$ est tel que $\kappa(C(Y)/\Delta(c_Y))=0$. 
(Voir [C01; \S 3] pour ces notions et r\'esultats).

Dans cette situation,  $S:=C(Y)$ est donc une surface K\" ahl\'erienne, et $\Delta(c_X):=D$ un diviseur orbifolde tel que 
$\kappa(S,K_S+D)=0$. En particulier: ou bien $\kappa(S)=0$, ou bien $S$ est rationnelle (le cas o\`u $S$ est 
r\'egl\'ee avec $q(S)>0$ est trait\'e dans [C01]). Nous allons traiter ici le cas le plus facile, o\`u $\kappa(S)=0$.

\begin{theorem} Soit $f:Y\ra S$  une fibration elliptique d'une
 vari\'et\'e K\" ahl\'erienne compacte $Y$ de dimension $3$  sur une 
surface (lisse) $S$. On suppose que $\kappa(S)=\kappa(S,K_S+D)=0$, o\`u $D:=\Delta(f)$ est 
le diviseur orbifolde des fibres multiples de $f$ (voir [C01]). Alors:

{\bf (1)} $\pi_1(X)$ est presque ab\'elien.

{\bf (2)} Si $S$ est projective, si $D$ est un diviseur exceptionnel de $S$, si l'orbifolde $S'$ obtenue par contraction de $D$ dans $S$ est bonne, et 
si $x,y,z\in Y$ \'etant trois points g\'en\'eriques, on peut trouver deux applications holomorphes $h_j:\Bbb C\ra Y, j=1,2$ 
telles que: $f(x),f(y)\in (f\circ h_1)(\Bbb C)$ et: $f(z),f(y)\in (f\circ h_2)(\Bbb C)$.

{\bf (3)} $d_Y\equiv 0$ si $S$ est projective, et si $S'$ est bonne.
\end{theorem}

{\bf D\'emonstration:} Dans ce cas, apr\`es rev\^etement \'etale fini, $S$ est bim\'eromorphe 
soit \`a une surface $K3$, soit \`a un tore complexe. 
Des arguments faciles montrent que $D$ est exceptionnel, puisque l'on a aussi: $\kappa(S,K_S+D)=0$.

La premi\`ere assertion r\'esulte alors de \ref{unifk3} et du lemme \ref{pi1ab} ci-dessous. La derni\`ere 
assertion est une cons\'equence imm\'ediate de la seconde et de la continuit\'e 
de $d_Y$ relativement \`a la topologie m\'etrique. 
Pour d\'emontrer la seconde propri\'et\'e, on consid\`ere un mod\`ele lisse du
 produit fibr\'e $Y':=Y\times_{T}$, si 
$T$ est une uniformisation lisse de l'orbifolde $S'$ obtenue en contractant le diviseur exceptionnel $D$ de $S$. 
On peut donc supposer (quitte \`a remplacer $Y$ par $Y'$) que $S=T$ est lisse, et est soit une surface $K3$, soit une
surface ab\'elienne. 

Dans le premier cas, on conclut en appliquant [B-L00] aux images r\'eciproques par $f$ 
des membres g\'en\'eriques de deux familles de 
courbes elliptiques de $S=T$ qui ne rencontrent pas $D$, la premi\`ere (resp. seconde) 
courbe elliptique contenant $x$ et $y$ (resp. $z$ et $y$). Dans le second cas, on applique [B-L00] au dessus (par $f$) de l'image 
d'une courbe enti\`ere $h:\Bbb C\ra T$ qui \'evite les points (en nombre fini) qui 
sont au-dessus des points singuliers de $S'$, et passe par $f(x), f(y)$ et $f(z)$ 
(une seule courbe enti\`ere suffit, dans ce cas) $\square$

\begin{lemma}\label{pi1ab} Soit $f: Y\ra S$ une fibration holomorphe 
entre vari\'et\'es complexes compactes et connexes, $Y$ K\" ahler. 
Si $F$ est la fibre g\'en\'erique de $f$, on suppose que:

{\bf (1)} $\pi_1(F)$ est presque ab\'elien.

{\bf (2)} Il existe $S^*\subset S$, un ouvert de Zariski dense tel que:

{\bf (2.a)} $\pi_1(S^*)$ est presque ab\'elien.

{\bf (2.b)} Pour tout $s\in S^*$, la fibre de $f$ au-dessus de $s$ est 
non-multiple (au sens des multiplicit\'es d\'efinies par $pgcd$; voir [C01, \S 9]).

Alors $\pi_1(Y)$ est presque nilpotent (ie: a un sous-groupe d'indice fini nilpotent), 
et presque ab\'elien si $Y$ est une vari\'et\'e sp\'eciale.
\end{lemma}

{\bf D\'emonstration:} Soit $Y^*:=f^{-1}(S^*)$. Le morphisme de 
groupes: $\pi_1(Y^*)\ra \pi_1(Y)$ d\'eduit de l'injection de $Y^*$ dans $Y$ est surjectif. Par l'hypoth\`ese 
(2.b), on a une suite exacte de groupes: 

$$1\ra \pi_1(F)\ra \pi_1(Y^*)\ra \pi_1(S^*)\ra 1,$$

d'o\`u l'on d\'eduit que $\pi_1(Y^*)$ est presque polycyclique (ie: a un sous-groupe d'indice fini polycyclique). 
Donc $\pi_1(Y)$ est presque polycyclique. De [A-N99], ou de [C01'], on d\'eduit que $\pi_1(Y)$ est presque nilpotent.

De [C01], on d\'eduit enfin que $\pi_1(Y)$ est presque ab\'elien si $Y$ est sp\'eciale $\square$
\bigskip

\section{\bf BIBLIOGRAPHIE}

[A03] F.Ambro. The Moduli B-Divisor of an LC-trivial Fibration. math.AG/0308143.

[A-N99] D.Arapura-M. Nori. Fundamental Groups of Algebraic Varieties and K\" ahler Manifolds. Comp. Math. 116 (1999), 173-193.

[Au78] T.Aubin. Equations de type Monge-Amp\`ere sur les vari\'et\'es K\" ahl\'eriennes compactes. Bull. Sc. Math. 102 (1978), 63-95.

[BPV 84] W.Barth-C.Peters-A.Van de Ven. Complex Surfaces. Erg. der Math. 4 Springer Verlag (1984)

[Bea 83] A.Beauville. Vari\'et\'es K\" ahl\'eriennes \`a premi\`ere classe de Chern nulle. J.Diff.Geom. 18 (1983), 755-782.

[De87] A.Besse. Einstein Manifolds. Erg. der Math. 10 (1987). Springer Verlag.

[Bog74] F.Bogomolov. The Decomposition of K\" ahler Manifolds with Trivial Canonical Class. Math. USSR Sb. 22 (1974), 580-583.

[Bog74] F.Bogomolov. Hamiltonian K\" ahler Manifolds. Sov. Math. Dokl. 19 (1978), 1462-1465.

[Bo93] J. Borzellino. Orbifolds of Maximal Diameter. Indiana Univ. Math. J. 42 (1993), 37-53.

[Bo-Z] J.Borzellino-S.H.Zhu. The Splitting Theorem for Orbifolds. Ill. J. Math. 38 (1994), 679-691.

[B-L00] J.Buzzard-S.Lu. Algebraic Surfaces Holomorphically Dominable by $\Bbb C^2$. Inv. Math. 139 (2000), 617-659.

[C 01] F.Campana. Special Manifolds and Classification Theory. Alg.geom. 0110051. 

Une introduction est expos\'ee dans: Special Varieties and Classification Theory: An Overview. 
Acta Applicandae Mathematicae. Kluwer Academic Publishers. 1-21. (2002)

[C01'] F.Campana. Ensembles de Green-Lazarsfeld et quotients r\'esolubles des groupes de K\" ahler. J.Alg. Geom. 10 (2001), 599-622.

[C-P02] F.Campana-T. Peternell. Projective Manifolds with Splitting Tangent Bundle, I. Math. Zeit. 241 (2002), 613-637.

[C-K-O 03] F.Catanese-J.H.Keum-K.Oguiso. Some remarks on the universal cover of an open $K3$ surface. Math. Ann. 325 (2003), 279-286.
 
[C-G71] J.Cheeger- Gromoll. The Splitting Theorem For Manifolds of Nonnegative Ricci Curvature. J.Diff. Geom. 6 (1971), 119-128.

[D-K00] J.P. Demailly-J. Koll\`ar. Semi-Continuity of Complex Singularity Exponents 
and K\" ahler-Einstein Metrics On Fano Orbifolds. arXiv.math. AG/9910118.

[F83] A. Fujiki. On Primitively Symplectic Compact K\"a hler V-Manifolds of Dimension Four. In: Classification of Algebraic and 
Analytic Manifolds. K. Ueno Ed. Progress in Math. 39 (1983), 71-250. Birkh\" auser Verlag.

[F-K-L] A.Fujiki-R.Kobayashi-S.Lu. On the Fundamental Group of Certain Open Normal Surfaces. Saitama Math. J. 11 (1993), 15-20.

[K-Z00] J.Keum-D.Q.Zhang. Algebraic Surfaces with Quotient Singularities. To appear in Proc. Symp. Geom. in East Asia. Aug. 2000.

[K-Z02] J. Keum and D. -Q. Zhang, Fundamental groups of open K3 surfaces, Enriques surfaces 
     and Fano 3-folds, Journal of Pure and Applied Algebra, 170 (2002), 67 -- 91.

[Li55] A. Lichnerowicz. Th\'eorie globale des connexions et des groupes d'holonomie (1955) Edizione Cremonese.
 Monografie matematiche 2. (Consiglio nazionale delle Ricerche)

[Mo-Mu82] S. Mori-S. Mukai. Uniruledness of the Moduli Space of Curves of Genus 11. LNM 1016 (1982), 334-353.

[P97] T.Peternell. Minimal Varieties with Trivial Canonical Classes, I. Math. Zeit. 217 (1997), 377-405.

[S-Z01] I. Shimada and D. -Q. Zhang, Classification of extremal elliptic K3 surfaces
       and fundamental groups of open K3 surfaces, Nagoya 
       Mathematical Journal, 161 (2001), 23--54. (Japan).

[Siu 87] Y.T.Siu. K\" ahler-Einstein Metrics. DMV Seminar 8. Birkh\" auser (1987).

[Y 78] S.T.Yau. On the Ricci Curvature of Compact K\" ahler Manifolds and the Complex Monge-Amp\`ere Equation. 
Comm. Pure and Applied Math. 31 (1978), 339-411.

\bigskip

F.Campana

D\'epartement de Math\'ematiques.

Universit\'e Nancy 1.

BP 239

F. 54506. Vandoeuvre-Les-Nancy. C\'edex.

e-mail: campana@iecn.u-nancy.fr

\end {document}